
\magnification1200
\input amstex.tex
\documentstyle{amsppt}

\vsize=18cm

\footline={\hss{\vbox to 2cm{\vfil\hbox{\rm\folio}}}\hss}
\nopagenumbers
\def\DJ{\leavevmode\setbox0=\hbox{D}\kern0pt\rlap
{\kern.04em\raise.188\ht0\hbox{-}}D}
\loadbold
\def\txt#1{{\textstyle{#1}}}
\baselineskip=13pt
\def\hf{{\textstyle{1\over2}}}
\def\a{\alpha}
\def\d{{\,\roman d}}
\def\e{\varepsilon}
\def\f{\varphi}
\def\G{\Gamma}
\def\k{\kappa}

\def\={\;=\;}

\def\zt{\zeta(\hf+it)}

\def\D{\Delta}
\def\no{\noindent}
 
\def\z{\zeta}

\def\hf{{\textstyle{1\over2}}}
\def\txt#1{{\textstyle{#1}}}
\def\f{\varphi}

\font\tenmsb=msbm10
\font\sevenmsb=msbm7
\font\fivemsb=msbm5
\newfam\msbfam
\textfont\msbfam=\tenmsb
\scriptfont\msbfam=\sevenmsb
\scriptscriptfont\msbfam=\fivemsb
\def\Bbb#1{{\fam\msbfam #1}}

\def \RR {\Bbb R}

\font\ff=cmr8
\def\txt#1{{\textstyle{#1}}}
\baselineskip=13pt

\font\teneufm=eufm10
\font\seveneufm=eufm7
\font\fiveeufm=eufm5
\newfam\eufmfam
\textfont\eufmfam=\teneufm
\scriptfont\eufmfam=\seveneufm
\scriptscriptfont\eufmfam=\fiveeufm
\def\mathfrak#1{{\fam\eufmfam\relax#1}}

\font\tenmsb=msbm10
\font\sevenmsb=msbm7
\font\fivemsb=msbm5
\newfam\msbfam
     \textfont\msbfam=\tenmsb
      \scriptfont\msbfam=\sevenmsb
      \scriptscriptfont\msbfam=\fivemsb
\def\Bbb#1{{\fam\msbfam #1}}

\def \RR {\Bbb R}

  \def\rightheadline{{\hfil{\ff
  On the Riemann zeta-function and the divisor problem II}\hfil\tenrm\folio}}

  \def\leftheadline{{\tenrm\folio\hfil{\ff
   A. Ivi\'c }\hfil}}
  \def\emptyheadline{\hfil}
  \headline{\ifnum\pageno=1 \emptyheadline\else
  \ifodd\pageno \rightheadline \else \leftheadline\fi\fi}

\topmatter
\title
ON THE RIEMANN ZETA-FUNCTION AND THE DIVISOR PROBLEM II
\endtitle
\author   Aleksandar Ivi\'c  \endauthor
\address
Aleksandar Ivi\'c, Katedra Matematike RGF-a
Universiteta u Beogradu, \DJ u\v sina 7, 11000 Beograd,
Serbia and Montenegro.
\endaddress
\email {\tt
ivic\@rgf.bg.ac.yu,  aivic\@matf.bg.ac.yu} \endemail
\dedicatory
Central European Journal of Mathematics 3(2) (2005), 203-214
\enddedicatory
\abstract
{Let $\D(x)$ denote the error term in the Dirichlet
divisor problem, and $E(T)$ the error term in the asymptotic
formula for the mean square of $|\zt|$. If
$E^*(t) = E(t) - 2\pi\D^*(t/2\pi)$ with $\D^*(x) =
 -\D(x)  + 2\D(2x) - \hf\D(4x)$, then we obtain
$$
\int_0^T |E^*(t)|^5\d t \;\ll_\e\; T^{2+\e}
$$
and
$$
\int_0^T |E^*(t)|^{544\over75}\d t \;\ll_\e\; T^{{601\over225}+\e}.
$$
It is also shown how bounds for moments of $|E^*(t)|$ lead to bounds
for moments of $|\zt|$.}
\endabstract
\endtopmatter

\head
1. Introduction and statement of results
\endhead

This work is the continuation of [8], where
several aspects of the connection between the divisor problem
and $\z(s)$, the zeta-function of Riemann, were investigated.
As usual, let
$$
\D(x) \;=\; \sum_{n\le x}d(n) - x(\log x + 2\gamma - 1)\eqno(1.1)
$$
denote the error term in the Dirichlet divisor problem, and
$$
E(T) \;=\;\int_0^T|\zt|^2\d t - T\left(\log\bigl({T\over2\pi}\bigr) + 2\gamma - 1
\right),\eqno(1.2)
$$
where $d(n)$ is the number of divisors of
$n, \gamma = -\G'(1) = 0.577215\ldots\,$
is Euler's constant. Instead of $\D(x)$ we work
with the modified function $\D^*(x)$ (see  M. Jutila [10]), where
$$
\D^*(x) \;:=\; -\D(x)  + 2\D(2x) - \hf\D(4x).\eqno(1.3)
$$
M. Jutila (op. cit.) investigated both the
local and global behaviour of the difference
$$
E^*(t) \;:=\; E(t) - 2\pi\D^*\bigl({t\over2\pi}\bigr),
$$
and in particular he proved that
$$
\int_0^T(E^*(t))^2\d t \;\ll\; T^{4/3}\log^3T.\eqno(1.4)
$$
In [8] this bound was complemented with the new bound
$$
\int_0^T (E^*(t))^4\d t \;\ll_\e\; T^{16/9+\e};\eqno(1.5)
$$
neither (1.4) or (1.5) seem to imply each other.
Here and later $\e$ denotes positive constants which are arbitrarily
small, but are not necessarily the same ones at each occurrence.
Our first aim is to obtain another bound for moments of $|E^*(t)|$.
This is given by

\bigskip
THEOREM 1. {\it We have}
$$
\int_0^T |E^*(t)|^5\d t \;\ll_\e\; T^{2+\e}.  \eqno(1.6)
$$

\bigskip
\no
From (1.4), (1.6) and H\"older's inequality for integrals, it follows that
$$
\eqalign{
\int_0^T |E^*(t)|^4\d t &= \int_0^T
|E^*(t)|^{2/3}|E^*(t)|^{10/3}\d t
\cr& \le \left(\int_0^T |E^*(t)|^2\d t\right)^{1/3}
\left(\int_0^T |E^*(t)|^5\d t\right)^{2/3}
\cr& \ll_\e T^{16/9+\e},\cr}
$$
which implies (1.5). This means that (1.6)  and
(1.4) together are stronger than  (1.5).
Another result of a more general nature (for the definition and
properties of exponent pairs see e.g., [3] or [6, Chapter 2]) is contained in

\bigskip
THEOREM 2. {\it Let $(\k, \lambda)$ be an exponent pair such that
$2\lambda \le 1+\k$, and
$$
V \;\ge \; T^{{1+\lambda-2\k\over3(2-\k)}+\e}.\eqno(1.7)
$$
Let $t_r\in [T,\,2T]\;(r = 1,\ldots, R)$ be points such that
$|t_r-t_s|\ge V\;(r\ne s)$ and $|E^*(t_r)| \ge V\;(r = 1,\ldots,R)$. Then
$$
R \ll_\e T^{1+\e}V^{-3} + T^{{1+4\k+\lambda\over3\k}+\e}V^{-{3\k+2\over\k}}.
\eqno(1.8)
$$}

\smallskip
From Theorem 2 we can obtain  specific bounds for moments of
$|E^*(t)|$, provided we choose the exponent pair $(\k,\lambda)$
appropriately. The optimal choice of the exponent pair is
hard to determine, since several conditions have to hold (see e.g., (5.5)).
However, by trying some of the standard exponent pairs one
can obtain a bound which is not far from the optimal bound
that the method allows. For instance, with the
exponent pair $(\k,\lambda) = (75/197, 104/197)$
(this exponent pair arises, in the terminology of exponent pairs,
as $\,(75/197, 104/197) = BA^3BA^3B(0,1)\,$) we can obtain

\bigskip
THEOREM 3.
{\it We have}
$$
\int_0^T |E^*(t)|^{544\over75}\d t \;\ll_\e\; T^{{601\over225}+\e}.
\eqno(1.9)
$$

\bigskip
One of the main reasons for investigating power moments of $|E^*(t)|$
is the possibility to use them to derive results on power moments
of $|\zt|$, which is one of the main themes in the theory of $\z(s)$.
A result in this direction is given by

\bigskip
THEOREM 4.
{\it Let $k \ge 1$ be a fixed real, and let $c(k)$ be such a
constant for which
$$
\int_0^T|E^*(t)|^{k}\d t \;\ll_\e\; T^{c(k)+\e}.\eqno(1.10)
$$
Then we have
$$
\int_0^T|\zt|^{2k+2}\d t \;\ll_\e\; T^{c(k)+\e}.\eqno(1.11)
$$}

\medskip

The constant $c(k)$ must satisfy
$$
c(k) \;\ge\; 1.\eqno(1.12)
$$
This is obvious if $k$ is an integer, as it follows from [6, Theorem 9.6].
If $k$ is not an integer, then this result yields ($p = {2k+2\over2[k]+2} > 1$)
$$
T \ll \int_0^T|\zt|^{2[k]+2}\d t \le \left(\int_0^T|\zt|^{2k+2}\d t\right)^{1/p}
T^{1-1/p}
$$
by H\"older's inequality for integrals. After simplification (1.12) easily follows again.

\medskip
{\bf Corollary 1.} We have
$$
\int_0^T|\zt|^{12}\d t \;\ll_\e \;T^{2+\e}.\eqno(1.13)
$$

\medskip\no
This follows from Theorem 1 and Theorem 4 (with $k=5$), and is the
well-known result of D.R. Heath-Brown [2], who had $\log^{17}T$ in
place of $T^\e$ on the right-hand side of (1.13).

\medskip
{\bf Corollary 2.} We have
$$
\int_0^T|\zt|^{1238\over75}\d t \;\ll_\e \;T^{{601\over225}+\e}.\eqno(1.14)
$$
\medskip\no
This follows from Theorem 3 and Theorem 4 (with
$k= {544\over75}$). The bound (1.14)
 does not follow from
(1.13) (and the strongest pointwise estimate for $|\zt|$), but on
the other hand (1.13) does not follow from (1.14). In principle, (1.14)
could be used for deriving zero-density bounds for $\z(s)$ (see e.g.,
[6, Chapter 10]), but very likely its use would lead to very
small improvements (if any) of the existing bounds.

\medskip
{\bf Acknowledgement}. I wish to thank Prof. Matti Jutila for valuable remarks.
\head
2. The necessary lemmas
\endhead

In this section we shall state the  lemmas which are necessary
for the proof of Theorem 1.

\medskip

LEMMA 1 (O. Robert--P. Sargos [11]). {\it Let $k\ge 2$ be a fixed
integer and $\delta > 0$ be given.
Then the number of integers $n_1,n_2,n_3,n_4$ such that
$N < n_1,n_2,n_3,n_4 \le 2N$ and}
$$
|n_1^{1/k} + n_2^{1/k} - n_3^{1/k} - n_4^{1/k}| < \delta N^{1/k}
$$
{\it is, for any given $\e>0$,}
$$
\ll_\e N^\e(N^4\delta + N^2).\eqno(2.1)
$$

This Lemma was crucial in obtaining the asymptotic formulas
for the third and fourth moment of $\D(x)$ in [9].
\medskip
LEMMA 2. {\it Let $T^\e \ll G \ll T/\log T$. Then we have}
$$
E^*(T) \le {2\over\sqrt{\pi}G}\int_0^{\infty} E^*(T+u)\,
{\roman e}^{-u^2/G^2}\d u + O_\e(GT^\e),\eqno(2.2)
$$
{\it and}
$$
E^*(T) \ge {2\over\sqrt{\pi}G}\int_0^\infty E^*(T-u)\,{\roman e}^{-u^2/G^2}\d u
+ O_\e(GT^\e).\eqno(2.3)
$$

\medskip \no
Lemma 2 follows on combining Lemma 2.2 and Lemma 2.3 of [8].

\medskip
The next lemma is F.V. Atkinson's classical explicit formula for $E(T)$ (see
[1], [6] or [7]).

\medskip
LEMMA 3. {\it Let $0 < A < A'$ be any two fixed constants
such that $AT < N < A'T$, and let $N' = N'(T) =
T/(2\pi) + N/2 - (N^2/4+ NT/(2\pi))^{1/2}$. Then }
$$
E(T) = \Sigma_1(T) + \Sigma_2(T) + O(\log^2T),\eqno(2.4)
$$
{\it where}
$$
\Sigma_1(T) = 2^{1/2}(T/(2\pi))^{1/4}\sum_{n\le N}(-1)^nd(n)n^{-3/4}
e(T,n)\cos(f(T,n)),\eqno(2.5)
$$
$$
\Sigma_2(T) = -2\sum_{n\le N'}d(n)n^{-1/2}(\log (T/(2\pi n))^{-1}
\cos(T\log (T/(2\pi n)) - T + \pi /4),\eqno(2.6)
$$
{\it with}
$$
\eqalign{\cr&
f(T,n) = 2T{\roman {arsinh}}\,\bigl(\sqrt{\pi n/(2T})\bigr) + \sqrt{2\pi nT
+ \pi^2n^2} - \pi/4\cr&
=  -\txt{1\over4}\pi + 2\sqrt{2\pi nT} +
\txt{1\over6}\sqrt{2\pi^3}n^{3/2}T^{-1/2} + a_5n^{5/2}T^{-3/2} +
a_7n^{7/2}T^{-5/2} + \ldots\,,\cr}\eqno(2.7)
$$
$$\eqalign{\cr
e(T,n) &= (1+\pi n/(2T))^{-1/4}{\Bigl\{(2T/\pi n)^{1/2}
{\roman {arsinh}}\,(\sqrt{\pi n/(2T})\,)\Bigr\}}^{-1}\cr&
= 1 + O(n/T)\qquad(1 \le n < T),
\cr}\eqno(2.8)
$$
{\it and $\,{\roman{arsinh}}\,x = \log(x + \sqrt{1+x^2}\,).$}

\medskip
LEMMA 4 (M. Jutila [10]). {\it If $A\in\RR$ is a constant, then we have}
$$
\cos\left(\sqrt{8\pi nT} +
\txt{1\over6}\sqrt{2\pi^3}n^{3/2}T^{-1/2} + A\right)
= \int_{-\infty}^\infty \a(u)\cos(\sqrt{8\pi n}(\sqrt{T} + u)
+ A)\d u,\eqno(2.9)
$$
{\it where $\a(u) \ll T^{1/6}$ for $u\not=0$,}
$$
\a(u) \ll T^{1/6}\exp(-bT^{1/4}|u|^{3/2}) \eqno(2.10)
$$
{\it for $u<0$, and}
$$
\a(u) = T^{1/8}u^{-1/4}\left(d\exp(ibT^{1/4}u^{3/2})
+ {\bar d}\exp(-ibT^{1/4}u^{3/2})\right) + O(T^{-1/8}u^{-7/4})\eqno(2.11)
$$
{\it for $u \ge T^{-1/6}$ and some constants $b\; (>0)$ and $d$.}

\head
3. The proof of  Theorem 1
\endhead

The proof is on the lines of [8]. We seek an upper bound
for $R$, the number of  points $\{t_r\}\in[T,2T]\,(r = 1,\ldots,\, R)$
such that $|E^*(t_r)| \ge V \ge T^\e$ and $|t_r-t_s|\ge V$ for
$r\ne s$. We consider separately the points where $E^*(t_r)$ is
positive or negative. Suppose the first case holds (the other one
is treated analogously). Then from Lemma 2 we have
$$ V \le E^*(t_r) \le
{2\over\sqrt{\pi}G}\int_0^\infty  E^*(t_r + u) \,{\roman
e}^{-u^2/G^2}\d u + O_\e(GT^\e).\eqno(3.1)
$$
The integral on the right-hand side is simplified by Atkinson's
formula (Lemma 3) and the truncated formula for $\D^*(x)$ (see
[8, eq. (6)]), as in [8].
We take $G = cVT^{-\e}$ (with sufficiently small $c>0$)
to make the $O$-term in (3.1) $\le \hf V$,
raise everything to the fourth power and sum over $r$.
By H\"older's inequality we obtain
$$
 RV^4 \ll_\e V^{-1}T^\e\max_{|u|\le G\log T}
\int_{T/2}^{2T} \f(t)\Bigl(\Sigma^4_4(X,N;u) +
\Sigma^4_5(X,N;u)+\Sigma^4_6(X;u)\Bigr)\d t,\eqno(3.2)
$$
with the notation introduced in (2.7), (2.8) and [8]:
$$\eqalign{\cr&
\Sigma_4(X,N;u) := t^{1/4}\sum_{X<n\le N}
(-1)^n d(n)n^{-3/4}e(t+u,n)\cos(f(t+u,n)),\cr&
\Sigma_5(X,N;u) := t^{1/4}
\sum_{X<n\le N}(-1)^n d(n)n^{-3/4}\cos(\sqrt{8\pi n(t+u)}-\pi/4),
\cr}\eqno(3.3)
$$
$$
\Sigma_6(X;u) := \sum_{n\le X}t^{-1/4}
(-1)^n d(n)n^{3/4} \cos(\sqrt{8\pi n(t+u)}-\pi/4).\eqno(3.4)
$$
Here we have $X = T^{1/3-\e},\, N = TG^{-2}\log T$, and
$\f(t)$ is a smooth, nonnegative function supported in
$\,[T/2,\,5T/2]\,$, such that $\f(t) = 1$ when $T \le t \le 2T$.
The basic idea is that the contributions of $\Sigma_6(X;u)$
and $\Sigma_5(X,N;u)$ will be approximately equal at $X$, and the
same will be true of $\Sigma_4(X,N;u)$ as well. In the latter case,
as was discussed in detail in [8], one has to use Lemma 4 to deal
with the complications arising from the presence of $\cos(f(t+u,n))$
in (3.3). The difference from [8] is that
the choice $G = cVT^{-\e}$ leads directly to (3.2), which
is in a certain sense optimal, while in [8] the choice
was $N = T^{5/9}$. Proceeding now as in [8] (here
Lemma 1 with $k=2$ was crucial) we obtain
$$
\eqalign{RV^4 &\ll_\e V^{-1}T^\e(T^{3/2}N^{1/2} + T^2X^{-1} +
T^{-1/2}X^{13/2} + X^5)\cr&
\ll_\e V^{-1}T^\e(T^2V^{-1} + T^{5/3})\cr&
\ll_\e T^{2+\e}V^{-2},\cr}\eqno(3.5)
$$
since $V < T^{1/3}$ in view of the best known estimates for $\D(x)$ and
$E(t)$. Namely with suitable $C>0$ one has (see M.N. Huxley [3], [4])
$$
\eqalign{
\D(x) &\;\ll\; x^{131/416}\log^Cx,\; 131/416 = 0.3149038\ldots,\cr
E(T) &\;\ll\; T^{72/227}\log^CT,\; 72/227 = 0.3171806\ldots.\cr}\eqno(3.6)
$$
Therefore (3.5) yields the large values estimate
$$
R \;\ll_\e \;T^{2+\e}V^{-6},
$$
and
Theorem 1 easily follows, as in [5] or   [6, Chapter 13] for moments of $\D(x)$.

\head
4. The proof of Theorem 2
\endhead

We start again from (3.1), choosing $G = cVT^{-\e} \,(< \hf V),\, T =t_r$,
so that we have
$$
 E^*(t_r) \ge V,\quad  E^*(t_r)\ll G^{-1}\int_0^\infty
 E^*(t_r+u){\roman e}^{-(u/G)^2}\d u
\eqno(4.1)
$$
in case $E^*(t_r) > 0$, and the case of negative values is analogous. We relabel the
points for which (4.1) holds in the sense that it will hold for $r = 1,\ldots, R$.
The proof is similar to the proof of (13.52) of Theorem 13.8 of [6]. To remove
the function $d(n)$ from the sums in
(3.3)--(3.4) we use the inequality (see the Appendix of [6])
$$
\sum_{r\le R}\,|({\boldsymbol \xi},{\boldsymbol \phi}_r)|^2 \;\le \;
||{\boldsymbol \xi}||^2\max_{r\le R}\sum_{s\le R}
\,|({\boldsymbol \phi}_r, {\boldsymbol\phi}_s)|,\eqno(4.2)
$$
where for two complex vector sequences ${\boldkey a} = \{a_n\}_{n=1}^\infty,
\,{\boldkey b} = \{b_n\}_{n=1}^\infty$ the inner product is defined as
$$
({\boldkey a},{\boldkey b}) \= \sum_{n=1}^\infty a_n{\bar b}_n.
$$
We shall also use (3.3)--(3.4) with $N = TG^{-2}\log T$. We shall consider separately
the points where $|\sum| \gg V$ when $\sum$ equals $\Sigma_4(X,N;u), \Sigma_5(X,N;u)$
or $\Sigma_6(X;u)\;(|u| \le G\log T)$, as the case may be. Taking the maximum over
$|u| \le G\log T$ over the whole sum, we may relabel the points such that they are
called again $t = t_r, \,r \le R$. Moreover, let $R_0$ denote the number of such
$t_r$'s (in each case) lying in an interval of length $T_0$, where $T_0$ is a function
of $V$ and $T$ that will be determined later. Thus $V\le T_0$ has to hold and
$$
R \;\ll\; R_0(1 + T/T_0).\eqno(4.3)
$$
As in the proof of Theorem 2, the choice of $X$ will be
$$
X \= T^{1/3-\e},
$$
when the largest term in $\sum_6$ is approximately equal to the smallest term in
$\sum_4$ and $\sum_5$. This choice exploits the specific structure of the function
$E^*(t)$, and leads to a better bound than was possible for large values of $\D(x)$
in Chapter 13 of [6]. Namely in the latter case the maximum occurred at $n = TG^{-2}\log T$,
but in our case $X = T^{1/3-\e} < TG^{-2}\log T$, since $V < T^{1/3-\e}$ must hold in view
of (3.6).
For example, from (3.4) and (4.2) (in case $|\sum_6| \gg V$ holds) we obtain
$$
\eqalign{
R_0V^2 &\ll {\log^2T\over\sqrt{T}} \max_{|u|\le G\log T, M\le X/2}
\sum_{r\le R_0}\left|\sum_{M<n\le2M}(-1)^nd(n)n^{3/4}{\roman e}^{i\sqrt{8\pi n(t_r+u})}
\right|^2\cr&
\ll {\log^2T\over\sqrt{T}}\max_{|u|\le G\log T,M\le X/2,r\le R_0}
M^{5/2}\log^4M\Biggl(M \,+ \cr&
\sum_{s\le R_0,s\ne r}\Bigl|\sum_{M<n\le2M}{\roman e}
^{i\sqrt{8\pi n}(\sqrt{t_r+u}-\sqrt{t_s+u})}\Bigr|\,\Biggr),\cr}\eqno(4.4)
$$
which corresponds to (13.60) of [6]. If we set
$$
f(x) = \sqrt{8\pi x}(\sqrt{t_r+u} - \sqrt{t_s+u}\,),
$$
then we can use the first derivative test (Lemma 2.1 of [6]) to deduce that the
contribution of $x = n$ (in the last sum in (4.4)) for which $|f'(x)| < 1/2$ is
$$
\eqalign{
& \ll \sum_{s\le R_0,s\ne r}{\sqrt{M}\over|\sqrt{t_r+u}-\sqrt{t_s+u}|}\cr&
\ll \sqrt{MT}\sum_{s\le R_0,s\ne r}{1\over|t_r-t_s|}\cr&
\ll \sqrt{MT}V^{-1}\log T,\cr}\eqno(4.5)
$$
since $|t_r-t_s| \ge V$ if $r\ne s$. The contribution of $|f'(x)| \ge 1/2$
is estimated by the theory of exponent pairs. The portion of the last sum in
(4.4) is, in this case,
$$
\ll R_0\left(|t_r-t_s|(MT)^{-1/2}\right)^{\k}M^\lambda \ll R_0T_0^\k M^{\lambda-\k/2}
T^{-\k/2},\eqno(4.6)
$$
since $|t_r-t_s|\le T_0$. Therefore from (4.4)--(4.6) it follows that
$$
\eqalign{
& R_0V^2 \ll T^{-1/2}X^{7/2}\log^6T + X^{3}V^{-1}\log^7T + R_0T_0^\k
X^{{5\over2}+\lambda
- {\k\over2}}T^{-{1\over2}-{\k\over2}}\log^6T\cr&
\ll T^{2/3}\log^6T + TV^{-1}\log^7T + R_0T_0^\k
T^{1+\lambda-2\k\over3}\log^6T.\cr}\eqno(4.7)
$$
The contribution of large values of $|\sum_4|$ and $|\sum_5|$ is estimated
analogously. We proceed, similarly as in (4.7), to obtain in these cases
$$
\eqalign{
R_0V^2 &\ll T^{1/2}\log^2T\max_{|u|\le G\log T,X\le M\le T^{1+\e}V^{-2}}\times\cr&\times
\sum_{r\le R_0}\left|\sum_{M<n\le2M}(-1)^nd(n)n^{-3/4}{\roman e}^{i\sqrt{8\pi n(t_r+u})}
\right|^2\cr&
\ll T^{1/2}\log^2T\max_{|u|\le G\log T,r\le R_0,X<M\le T^{1+\e}V^{-2}}
M^{-1/2}\log^4M\Biggl(M \,+ \cr&
\sum_{s\le R_0,s\ne r}\Bigl|\sum_{M<n\le2M}{\roman e}
^{i\sqrt{8\pi n}(\sqrt{t_r+u}-\sqrt{t_s+u})}\Bigr|\,\Biggr)
\cr&
\ll_\e T^{1+\e}V^{-1} + R_0T_0^\k T^{1/2-\k/2}\log^5T\max_{X<M\le T^{1+\e}V^{-2}}
M^{\lambda-\k/2-1/2}
.\cr}\eqno(4.8)
$$
The hypothesis in the formulation of the theorem was that
$$
2\lambda \;\le\; \k + 1,\;\eqno(4.9)
$$
hence by combining (4.7) and (4.8) it follows that
$$
R_0V^2 \ll_\e T^{1+\e}V^{-1} + R_0T_0^\k T^{{1+\lambda-2\k\over3}+\e},\eqno(4.10)
$$
since $T^{2/3} \le TV^{-1}$ because $V \le T^{1/3}$ has to hold. If we choose
$$
T_0 \= V^{2\over\k}T^{{2\k-1-\lambda\over3\k}-{2\e\over\k}}\eqno(4.11)
$$
then (4.10) reduces to $R_0V^2 \ll_\e T^{1+\e}V^{-1}$, and the condition
$T_0\ge V$ becomes
$$
V \;\ge\; T^{{1+\lambda-2\k\over3(2-\k)}+\e},\eqno(4.12)
$$
which is (1.7). Therefore (4.10) gives
$$
R \ll R_0(1 + T/T_0) \ll_\e T^{1+\e}V^{-3} + T^{{1+4\k+\lambda\over3\k}+\e}
V^{-{3\k+2\over\k}},
$$
thereby completing the proof of Theorem 2.

\medskip
\head
5. The proof of Theorem 3
\endhead
With the choice $(\k,\lambda) = (75/197,\,104/197)$ it is seen that
(1.7) and (1.8) of Theorem 2 reduce to
$$
R \ll_\e T^\e(TV^{-3} + T^{601\over225}V^{-{619\over75}})\quad
(V \ge T^{{151\over957}+\e}, \; {\txt{151\over957}} = 0.157784\ldots\,).\eqno(5.1)
$$
Let
$$
J_V(T)  \:= \Bigl\{\,t\in\,[T,\,2T]\;:\; V \le |E^*(t)| < 2V\Bigr\},
$$
and write
$$
\int_T^{2T}|E^*(t)|^{544\over75}\d t \;\ll_\e\; \log T\max_{V\ge T^\e}
\int_{J_V(T)}|E^*(t)|^{544\over75}\d t  + T^{1+\e}.\eqno(5.2)
$$
For $V \le T^{151/957+\e}$ we have, on using (1.6) of Theorem 1,
$$
\eqalign{&
\int_{J_V(T)}|E^*(t)|^{544\over75}\d t = \int_{J_V(T)}|E^*(t)|^{5}
|E^*(t)|^{169\over75}\d t\cr&
\ll_\e T^{2+{169\over75}\cdot{151\over957}+\e} \le T^{601\over225}.\cr}
\eqno(5.3)
$$
Suppose now that $V \ge T^{151/957+\e}$, and divide $\,[T,\,2T]$ into
subintervals of length $V$ (the last of these subintervals may be shorter).
Let $|E^*(\tau_j)|$ be the supremum of $|E^*(t)|$ in the $j$th of these
subintervals, and let further $t_1,...,t_{R_V}$ denote the $\tau_j$'s with
even or odd indices such that the intersection of the $j$th subinterval
and $J_V(T)$ is non-empty. Then $|t_r-t_s| \ge V$ for $r\ne s$, and (5.1) gives
$$
R_V \ll_\e T^\e(TV^{-3} + T^{601\over225}V^{-{619\over75}})
\ll_\e T^{{601\over225}+\e}V^{-{619\over75}}\eqno(5.4)
$$
for
$$
V\;\le\; T^{188\over591},\quad  188/591 = 0.3181049\ldots\,.\eqno(5.5)
$$
But in view of (3.6) it is seen that (5.5) is always satisfied
(the choice of our exponent pair was made to ensure that this is
indeed the case), and we obtain from (5.4)
$$
\int_{J_V(T)}|E^*(t)|^{544\over75}\d t \ll R_VV^{1+{544\over75}}
\ll_\e T^{{601\over225}+\e}V^{-{619\over75}}V^{{619\over75}} =
T^{{601\over225}+\e}.
\eqno(5.6)
$$
Theorem 3 follows then from (5.2), (5.3) and (5.6), on replacing
$T$ by $T2^{-j}$ and summing over $j = 1,2,\dots\,$.

\medskip
\head
6. The proof of Theorem 4
\endhead
To prove Theorem 4 it is enough to prove that
$$
R \ll_\e T^{c(k)+\e}V^{-2k-2},\eqno(6.1)
$$
where $R$ is the number of points $t_r \in [T,\,2T]\,(r = 1,\ldots\,,R),$
such that $|\z(\hf+it_r)| \ge V$ with $|t_r-t_s| \ge 1$ for $r\ne s$ and
$V\ge T^\e$. We denote actually by $R$ the number of points with even
and odd indices, so that the intervals $[t_r - {1\over3}, t_r + {1\over3}]$
are disjoint.
Then we have, using Theorem 1.2 of [7] with $k=2, \delta = {1\over3}$,
$$
\eqalign{
RV^2 & \le \sum_{r=1}^R |\z(\hf+it_r)|^2
\ll \log T\sum_{r=1}^R\int_{t_r-{1\over3}}
^{t_r+{1\over3}}|\zt|^2\d t\cr&
\ll \log T\sum_{j=1}^J \int_{\tau_j-G}^{\tau_j+G}|\zt|^2\d t,\cr}\eqno(6.2)
$$
where $\tau_j \in [T-G,\,T+G]\;(j = 1,\ldots\,,J)$ is a system of points such that
$|\tau_j - \tau_\ell| \ge 2G$ for $j\ne \ell$ and $T^\e \le G = G(T) \ll T$.
By the definition of $E^*(t)$ we have
$$
\eqalign{&\int_{\tau_j-G}^{\tau_j+G}|\zt|^2\d t
= E(\tau_j+G) - E(\tau_j-G) + O(G\log T)\cr&
= E^*(\tau_j+G) - E^*(\tau_j-G) + 2\pi\D^*\left({\tau_j+G\over2\pi}\right)
- 2\pi\D^*\left({\tau_j-G\over2\pi}\right) + O(G\log T)\cr&
= E^*(\tau_j+G) - E^*(\tau_j-G) + O_\e(GT^\e).\cr}
$$
Here we used the fact that
$$
\D^*(x) - \D^*(y) \ll_\e x^\e(x - y + 1)\qquad(1 \ll y \le x),
$$
which follows from (1.1), (1.3) and $d(n) \ll_\e n^\e$.
 This arithmetic property of $d(n)$ is essential,
since it makes it possible to connect the large values of $|\zt|$ to sums
of values of $E^*(t)$, and hence to exploit the special structure of
the function $E^*(t)$. If we worked only with $E(t)$, we would obtain Theorem 4,
where (1.10) has $E^*(t)$ replaced by $E(t)$. However,
the existing estimates for the moments
of $|E(t)|$ (see  [5] and Chapter 13 of [6]) are not as strong as the moments
of $|E^*(t)|$ (cf. (1.6) and (1.9)).

Returning to the proof, note that (6.2) yields
$$
RV^2 \ll_\e \log T\left\{\sum_{j=1}^J(E^*(\tau_j+G) - E^*(\tau_j-G))\right\}
+ RGT^\e,
$$
giving
$$
RV^2 \ll_\e \log T\left\{\sum_{j=1}^J(E^*(\tau_j+G) - E^*(\tau_j-G))\right\}
\eqno(6.3)
$$
with
$$
G \= cV^2T^{-\e}\eqno(6.4)
$$
and sufficiently small $c>0$. If we use Lemma 2
we may replace $\sum_j E^*(\tau_j+G)$ by its majorant
$$
{2\over\sqrt{\pi}G}\int_0^\infty \sum_{j=1}^J E^*(\tau_j+G+u)
\,{\roman e}^{-u^2/G^2}\d u + RGT^\e,
$$
and similarly for the sum with $E^*(\tau_j-G)$. By H\"older's
inequality we have (since $J \le R$)
$$\eqalign{&
\int_0^\infty \sum_{j=1}^J E^*(\tau_j+G+u)\,{\roman e}^{-u^2/G^2}\d u\cr&
\ll \int_0^\infty \,{\roman e}^{-u^2/G^2}
\left(\sum_{j=1}^J|E^*(\tau_j+G+u)|^{k}\right)^{1\over k}
R^{1-{1\over k}}\d u\cr&
\ll R^{1-{1\over k}}\left(\int_0^\infty \,{\roman e}^{-u^2/G^2}
\sum_{j=1}^J|E^*(\tau_j+G+u)|^{k}\d u\right)^{1\over k}G^{1-{1\over k}}\cr&
\ll (GR)^{1-{1\over k}}\left(\int_{T/2}^{5T/2}|E^*(t)|^{k}\d t
\right)^{1\over k}\cdot\log T,\cr}\eqno(6.5)
$$
by breaking the system of points $\tau_j$ into $\ll\log T$ subsystems
with $|\tau_j-\tau_\ell| \ge G\log T$ for $\ell\ne j$.
From (1.10) and (6.3)--(6.5) it follows that
$$
RV^2 \ll_\e T^\e(RV^2)^{1-{1\over k}}\cdot T^{{c(k)\over k}+\e}V^{-2},
$$
which on simplifying yields
$$
R \ll_\e T^{c(k)+\e}V^{-2k-2},\eqno(6.6)
$$
and (6.6) implies easily (1.11) of Theorem 4. By the same method one also
obtains $\gamma \le c(k)/(k+1)$ for every $k \ge 1$, if
$$
\gamma :\= \inf\{\,g\ge 0\,:\, E^*(T) \ll T^g\,\},
$$
but better bounds for $\gamma$ can be derived from short
interval results on $E^*(t)$, provided they can be obtained.
The existing results make it hard to even conjecture what
should be the true value of $\gamma$.

\vfill
\eject
\topglue2cm
\Refs
\bigskip

\item{[1]} F.V. Atkinson, The mean value of the Riemann zeta-function,
Acta Math. {\bf81}(1949), 353-376.

\item{[2]} D.R. Heath-Brown, The twelfth power moment of the Riemann
zeta-function, Quart. J. Math. (Oxford) {\bf29}(1978), 443-462.

\item{[3]} M.N. Huxley, Area, Lattice Points and Exponential
Sums, Oxford Science Publications, Clarendon Press,
Oxford, 1996.

\item{[4]} M.N. Huxley,  Exponential sums and lattice points III, Proc.
London Math. Soc. 3{\bf87}(2003), 591-609.

\item{[5]} A. Ivi\'c, Large values of the error term in the
divisor problem, Invent. Math. {\bf71}(1983), 513-520.

\item{[6]} A. Ivi\'c, The Riemann zeta-function, John Wiley \&
Sons, New York, 1985 (2nd ed. Dover, Mineola, New York, 2003).

\item{[7]} A. Ivi\'c, The mean values of the Riemann zeta-function,
LNs {\bf 82}, Tata Inst. of Fundamental Research, Bombay (distr. by
Springer Verlag, Berlin etc.), 1991.

\item{[8]} A. Ivi\'c, On the Riemann zeta-function and the divisor problem,
Central European J. Math. {\bf(2)(4)} (2004), 1-15.

\item{[9]} A. Ivi\'c and P. Sargos, On the higher moments of the
error term in the divisor problem, to appear.

\item{[10]} M. Jutila, Riemann's zeta-function and the divisor problem,
Arkiv Mat. {\bf21}(1983), 75-96 and II, ibid. {\bf31}(1993), 61-70.

\item{[11]} O. Robert and P. Sargos, Three-dimensional
exponential sums with monomials, J. reine angew. Math. (in print).

\endRefs

\vskip3cm\sevenrm
Aleksandar Ivi\'c

Katedra Matematike RGF-a

Universitet u Beogradu, \DJ u\v sina 7

 11000 Beograd, Serbia and Montenegro

{\sevenbf ivic\@rgf.bg.ac.yu}

\vfill


\bye